\documentclass[a4paper,12pt,twoside]{article}
\usepackage[francais]{babel}
\usepackage{amssymb, amsmath, eucal}
\usepackage[all]{xy}
\usepackage{mathrsfs}
\begin{document}
\textwidth15.5cm
\textheight22.5cm
\voffset=-13mm
\newtheorem{The}{Theorem}[section]
\newtheorem{Lem}[The]{Lemma}
\newtheorem{Prop}[The]{Proposition}
\newtheorem{Cor}[The]{Corollary}
\newtheorem{Rem}[The]{Remark}
\newtheorem{Obs}[The]{Observation}
\newtheorem{Titre}[The]{\!\!\!\! }
\newtheorem{Conj}[The]{Conjecture}
\newtheorem{Question}[The]{Question}
\newtheorem{Prob}[The]{Problem}
\newtheorem{Def}[The]{Definition}
\newtheorem{Not}[The]{Notation}
\newtheorem{Claim}[The]{Claim}
\newtheorem{Conc}[The]{Conclusion}
\newtheorem{Ex}[The]{Example}
\newcommand{\C}{\mathbb{C}}
\newcommand{\R}{\mathbb{R}}
\newcommand{\N}{\mathbb{N}}
\newcommand{\Z}{\mathbb{Z}}
\newcommand{\Q}{\mathbb{Q}}

\begin{center}

{\Large\bf Limits of Moishezon Manifolds under Holomorphic Deformations}

\end{center}

\begin{center}

 {\large Dan Popovici}

\end{center}

\vspace{1ex}

\noindent {\small {\bf Abstract.} Given a (smooth) complex analytic family of compact complex manifolds, we prove that the central fibre must be Moishezon if the other fibres are Moishezon. Using a {\it strongly Gauduchon metric} on the central fibre whose existence was proved in our previous work on limits of projective manifolds, we show that the irreducible components of the relative Barlet space of divisors contained in the fibres are proper over the base even under the weaker assumption that the $\partial\bar\partial$-lemma hold on all the fibres except, possibly, the central one. This implies that the algebraic dimension of the central fibre cannot be lower than that of the generic fibre. Since the latter is already maximal thanks to the Moishezon assumption, the central fibre must be of maximal algebraic dimension, hence Moishezon.}

\vspace{3ex}

\section{Introduction}\label{section:introd}

 Let $\pi:{\cal X}\rightarrow\Delta$ be a complex analytic (also termed holomorphic) family of compact complex manifolds (in the sense of [Kod86]) over a ball $\Delta$ about the origin in some $\C^m$. This means that ${\cal X}$ is a complex manifold and $\pi$ is a proper holomorphic submersion. Thus $\pi$ is a smooth map in the sense of algebraic geometry and the fibres $X_t=\pi^{-1}(t)$, $t\in\Delta$, are (smooth) compact complex manifolds whose common complex dimension will be denoted by $n$. On the other hand, recall that a compact complex manifold $X$ is said to be {\it Moishezon} if it admits a holomorphic {\it modification} (i. e. a holomorphic bimeromorphic map) $\mu : \widetilde{X}\rightarrow X$ from some projective manifold $\widetilde{X}$ (cf. [Moi67]). Bringing these two notions together, we set out to prove the following statement.

\begin{The}\label{The:moidef} Let $\pi:{\cal X}\rightarrow\Delta$ be a complex analytic family of compact complex manifolds such that $X_t$ is Moishezon for every $t\in\Delta^{\star}:=\Delta\setminus\{0\}$. Then $X_0$ is again Moishezon.

\end{The}

 This result generalises our main result of [Pop09] where the same conclusion was obtained under the stronger assumption that $X_t$ be projective for every $t\in\Delta^{\star}$. However, the proof of the present Theorem \ref{The:moidef} will make crucial use of a special class of Gauduchon metrics that were introduced under the name of {\it strongly Gauduchon metrics} in [Pop09] and already played a key role there. Such a metric was shown to exist on $X_0$ if the $\partial\bar{\partial}$-lemma was assumed to hold on all the other fibres $X_t$ with $t\in\Delta^{\star}$ (cf. [Pop09, Proposition 4.1]). One can even find a family $(\gamma_t)_{t\in\Delta}$, varying in a $C^{\infty}$ way with $t$, of {\it strongly Gauduchon metrics} on the fibres $(X_t)_{t\in\Delta}$. Such a family enables one to uniformly bound the volumes of the divisors that form an arbitrary irreducible component of the relative Barlet space ${\cal C}^{n-1}({\cal X}/\Delta)$ of effective analytic divisors $Z_t$ contained in the fibres $X_t$ over any compact subset of $\Delta$. It follows that the irreducible components of ${\cal C}^{n-1}({\cal X}/\Delta)$ are proper over $\Delta$ in the following sense.

\begin{Prop}\label{Prop:reldivprop} Let $\pi:{\cal X}\rightarrow\Delta$ be a complex analytic family of compact complex manifolds such that the $\partial\bar{\partial}$-lemma holds on $X_t$ for every $t\in\Delta^{\star}$. Then the canonical holomorphic projection

$$\mu_{n-1}:{\cal C}^{n-1}({\cal X}/\Delta)\rightarrow\Delta, \hspace{3ex} \mu_{n-1}(Z_t)=t,$$

\noindent mapping every divisor $Z_t\subset X_t$ contained in some fibre $X_t$ to the base point $t\in\Delta$, has the property that its restrictions to the irreducible components of ${\cal C}^{n-1}({\cal X}/\Delta)$ are proper.

 \end{Prop}

  While the irreducible components of the Barlet space of cycles of arbitrary codimension ${\cal C}(X)$ need not be compact on a general compact complex manifold $X$ (cf. [Lie78]), compactness of the irreducible components of the Barlet space ${\cal C}^{n-1}(X)$ of divisors of $X$ always holds if $X$ is compact (see e.g. [CP94, Remark 2.18.]). Thus the absolute case of Proposition \ref{Prop:reldivprop} (i.e. when $\Delta$ is reduced to a point) is well-known and no special assumption is necessary. However, the relative counterpart fails in general as shown by an example given by Fujiki and Pontecorvo [FP09] of a family of compact non-K\"ahler complex surfaces of class VII in which the algebraic dimension drops from $1$ to $0$ on the central fibre. In particular, properness cannot hold for the irreducible components of the relative space of divisors.\footnote{The author is grateful to Fr\'ed\'eric Campana for pointing out to him this example of [FP09].} It is thus owing to the $\partial\bar{\partial}$-lemma assumption on the fibres above $\Delta^{\star}$ that Proposition \ref{Prop:reldivprop} holds. Notice that, since the only compact complex surfaces on which the $\partial\bar{\partial}$-lemma holds are the K\"ahler ones, the family exhibited in [FP09] does not satisfy the hypothesis of Proposition \ref{Prop:reldivprop}.

 Here is how Theorem \ref{The:moidef} follows from Proposition \ref{Prop:reldivprop}. The latter certainly applies to the family considered in the former since the $\partial\bar\partial$-lemma holds on every Moishezon manifold. Properness guarantees that the images of the irreducible components of ${\cal C}^{n-1}(X/\Delta)$ under $\mu_{n-1}$ are analytic subsets of $\Delta$ thanks to Remmert's Proper Mapping Theorem. Let $\Sigma_{\nu}\subsetneq\Delta$, for $\nu\in\Z$, be those such images (at most countably many) that are {\it strictly} contained in $\Delta$. Each $\Sigma_{\nu}$ is thus a proper analytic subset of $\Delta$. Bearing in mind the structure of the irreducible components of the (relative) Barlet space of cycles as described in [Bar75], we see that every irreducible component $S$ of ${\cal C}^{n-1}({\cal X}/\Delta)$ gives rise to an analytic family (in the sense of [Bar75, Th\'eor\`eme 1, p. 38]) of relative effective divisors $(Z_s)_{s\in S}$ such that $Z_s\subset X_{\mu_{n-1}(s)}$ for all $s\in S$. We can either have

\begin{equation}\label{eqn:fam-cycles1}\mu_{n-1}(S)=\Delta \hspace{3ex} \mbox{or}\end{equation}


\begin{equation}\label{eqn:fam-cycles2}\mu_{n-1}(S)=\Sigma_\nu\subsetneq\Delta, \hspace{2ex} \mbox{for some}\,\,\, \nu\in\Z.\end{equation}

 Let $\Sigma=\bigcup_{\nu}\Sigma_{\nu}\subsetneq\Delta$. Thus every divisor $Z_{s_0}$ contained in a fibre $X_{t_0}$ lying above some point $t_0=\mu_{n-1}(s_0)\in\Delta\setminus\Sigma$ (call such a fibre {\it generic}) stands in an analytic family of divisors $(Z_s)_{s\in S}$ covering the whole base $\Delta$ as in (\ref{eqn:fam-cycles1}) (call these divisors {\it generic}), while the {\it exceptional} fibres $X_t$ (i.e. those above points $t\in\Sigma$) may have extra divisors (those standing in {\it isolated} families satisfying (\ref{eqn:fam-cycles2})) besides the {\it generic} divisors that ``sweep'' $\Delta$ in families with the property (\ref{eqn:fam-cycles1}).  

 In other words, properness of the irreducible components of ${\cal C}^{n-1}({\cal X}/\Delta)$ ensures that every fibre (in particular $X_0$) has at least as many divisors (the {\it generic} ones) as the {\it generic} fibres of the family. 
On the other hand, the algebraic dimension of any fibre $X_t$ is the maximal number of effective prime divisors meeting transversally at a generic point of $X_t$ (see e.g. [CP94, Remark 2.22]). It follows from the last two assertions that the algebraic dimension of $X_0$ is $\geq$ the algebraic dimension of the {\it generic} fibre. However, the algebraic dimension of any $X_t$ with $t\neq 0$ is maximal (i.e. equals the complex dimension $n$) thanks to the Moishezon assumption (known to be equivalent to the maximality of the algebraic dimension by [Moi67]). Thus the algebraic dimension of $X_0$ must be maximal or, equivalently, $X_0$ must be Moishezon. 

 The analytic cycle approach adopted in the present work offers an alternative to the K\"ahler metric approach of [Pop09] when the undertaking is aimed at proving that the limit fibre $X_0$ is Moishezon. However, the method of [Pop09] relying on {\it singular Morse inequalities} will most likely prove vital in a future attack on the standard conjecture predicting that the deformation limit of a holomorphic family of compact K\"ahler (or merely {\it class} ${\cal C}$) manifolds is {\it class} ${\cal C}$. As explained in the introduction to [Pop09], the only missing link in this direction is Demailly's conjecture on {\it transcendental Morse inequalities}.

\section{Proof of Proposition \ref{Prop:reldivprop}}\label{section:reldivprop}

 To put the result stated in Proposition \ref{Prop:reldivprop} in context, we feel bound to make a few comments. Recall that a compact complex manifold $X$ is said to be in the {\it class} ${\cal C}$ if it admits a holomorphic modification $\mu : \widetilde{X}\rightarrow X$ from a compact K\"ahler manifold $\widetilde{X}$ (cf. e.g. [Dem97, chapter VI, $\S. 12$]). {\it Class} ${\cal C}$ {\it manifolds} were introduced by Fujiki in [Fuj78] as meromorphic images of compact K\"ahler manifolds; they were subsequently given the nice description adopted as a definition above by Varouchas in [Var86]. It has been known since the work of Fujiki (see [Fuj78, Theorem 4.9.]) that the irreducible components of the Barlet space of cycles ${\cal C}(X)$ of a {\it class} ${\cal C}$ {\it manifold} $X$ are compact. (They are even {\it class} ${\cal C}$ by [Cam80, Corollaire 3], but this extra property is immaterial to our purposes here.) As already mentioned, this last property fails if $X$ is merely supposed to be compact (although it holds for divisors), while the {\it class} ${\cal C}$ assumption is the minimal requirement on $X$ that we are aware of ensuring compactness of the irreducible components. 

 It thus appears natural to conjecture the (more general) relative case. 

\begin{Conj}\label{Conj:relpropC} Let $\pi:{\cal X}\rightarrow\Delta$ be a complex analytic family of compact complex manifolds such that the fibre $X_t:=\pi^{-1}(t)$ is a {\it class} ${\cal C}$ manifold for every $t\in\Delta$. Then the irreducible components of the relative Barlet space ${\cal C}({\cal X}/\Delta)$ of cycles on ${\cal X}$ are proper over $\Delta$.

\end{Conj}

 We have used the standard notation

 $${\cal C}({\cal X}/\Delta)=\bigcup\limits_{0\leq p\leq n}{\cal C}^p({\cal X}/\Delta),$$

\noindent where ${\cal C}^p({\cal X}/\Delta)$ stands for the relative Barlet space of effective analytic $p$-cycles contained in the fibres $X_t$. The special case of the above conjecture where all the fibres are supposed to be K\"ahler is well-known and quite easy to prove, but the general case of {\it class} ${\cal C}$ fibres is still elusive. We may even go so far as conjecture the same conclusion when the {\it class} ${\cal C}$ assumption is made to skip one of the fibres.

\begin{Conj}\label{Conj:relpropC*} Let $\pi:{\cal X}\rightarrow\Delta$ be a complex analytic family of compact complex manifolds such that the fibre $X_t:=\pi^{-1}(t)$ is a {\it class} ${\cal C}$ manifold for every $t\in\Delta^{\star}$. Then the irreducible components of the relative Barlet space of cycles ${\cal C}({\cal X}/\Delta)$ are proper over $\Delta$.

\end{Conj}

 Our Proposition \ref{Prop:reldivprop} answers affirmatively the stronger Conjecture \ref{Conj:relpropC*} in the special case of divisors (and even under the weaker $\partial\bar\partial$-lemma assumption which is known to hold on any {\it class} ${\cal C}$ manifold). A tantalising special case of Conjecture \ref{Conj:relpropC*} is the one where the fibres $X_t$ with $t\neq 0$ are supposed to be even K\"ahler. The central fibre $X_0$ is then expected to be {\it class} ${\cal C}$, but proving the compactness of the irreducible components of its Barlet space of cycles would be a first step towards confirming this expectation. 

 We will now outline the first moves towards possible solutions of these conjectures that will make the (considerable) difficulties apparent while proving Proposition \ref{Prop:reldivprop} by a crucial application of a result from [Pop09].  

 Fix a complex analytic family of compact complex manifolds $\pi:{\cal X}\rightarrow\Delta$ and let $n$ denote the complex dimension of the fibres $X_t$, $t\in\Delta$. Recall that all the fibres $X_t$, $t\in\Delta$, are {\it a fortiori} $C^{\infty}$-diffeomorphic to a fixed compact $C^{\infty}$-manifold $X$ and that only the complex structure $J_t$ of $X_t$ varies (holomorphically) with $t\in\Delta$ (see e.g. [Kod86]). Thus the De Rham cohomology groups $H^k_{DR}(X_t, \, \C)$ of the fibres can be identified with a fixed $H^k(X, \, \C)$ for all $t\in\Delta$, while the Dolbeault cohomology groups $H^{p, \, q}(X_t, \, \C)$ vary with the complex structure $J_t$. 

 For every $p\in\{0, 1, \dots , n\}$, consider the relative Barlet space ${\cal C}^p({\cal X}/\Delta)$ of effective analytic $p$-cycles on ${\cal X}$ that are contained in the fibres $X_t$. It is a subspace of the (absolute) Barlet space ${\cal C}^p({\cal X})$ of compact $p$-cycles on ${\cal X}$. Further recall that ${\cal C}({\cal X}):=\cup_p{\cal C}^p({\cal X})$ is the Chow scheme of ${\cal X}$ (which, by definition, parametrises the compactly supported analytic cycles of ${\cal X}$) that Barlet endowed with a natural structure as a Banach analytic set whose irreducible components are finite-dimensional analytic sets (cf. [Bar75]). Moreover, any irreducible component $S$ of ${\cal C}({\cal X})$ arises as an analytic family of compact cycles $(Z_s)_{s\in S}$ parametrised by $S$, while giving an analytic family $(Z_s)_{s\in S}$ of compact cycles of dimension $p$ on ${\cal X}$ is equivalent to giving an analytic subset

$${\cal Z}=\{(s,\, z)\in S\times {\cal X}\, /\, z\in |Z_s|\}\subset S\times{\cal X},$$

\noindent where $|Z_s|$ denotes the support of the cycle $Z_s$, such that the restriction to ${\cal Z}$ of the natural projection on $S$ is proper, surjective and has fibres of pure dimension $p$ (cf. [Bar75, Th\'eor\`eme 1, p. 38]). Recall finally Lieberman's strengthened form ([Lie78, Theorem 1.1]) of Bishop's Theorem [Bis64]: a subset $S\subset{\cal C}({\cal X})$ is relatively compact if and only if the supports $|Z_s|$, $s\in S$, all lie in a same compact subset of ${\cal X}$ and the $\widetilde{\omega}$-volume of $Z_s$ is uniformly bounded when $s\in S$ for some (hence any) Hermitian metric $\widetilde{\omega}$ on ${\cal X}$. Here, as usual, the $\widetilde{\omega}$-volume of a $p$-cycle $Z_s\subset{\cal X}$ is defined to be

$$v_{\widetilde{\omega}}(Z_s):=\int\limits_{{\cal X}}[Z_s]\wedge\widetilde{\omega}^p=\int\limits_{Z_s}\widetilde{\omega}^p,$$

\noindent where $[Z_s]$ is the current of integration on the cycle $Z_s$.

 Let us now fix $p\in\{0, 1, \dots , n\}$ and suppose that the family ${\cal X}=(X_t)_{t\in\Delta}$ satisfies the hypothesis of Proposition \ref{Prop:reldivprop}. So $X_t$ is merely assumed to have the $\partial\bar\partial$-lemma hold for every $t\in\Delta^{\star}$, while no special assumption is made on $X_0$. Fix also a family $(\gamma_t)_{t\in\Delta}$ of Hermitian metrics, varying in a $C^{\infty}$ way with $t$, on the respective fibres $(X_t)_{t\in\Delta}$. Let $(Z_t)_{t\in\Delta^{\star}}$ be a differentiable family of effective analytic $(n-p)$-cycles such that $Z_t\subset X_t$ for every $t\in\Delta^{\star}$. The main difficulty in proving the properness predicted by Conjecture \ref{Conj:relpropC*} is to ensure the uniform boundedness of the $\gamma_t$-volumes of the cycles $Z_t:$

$$v_{\gamma_t}(Z_t)=\int\limits_X[Z_t]\wedge\gamma_t^{n-p}, \hspace{3ex} t\in\Delta^{\star},$$

\noindent as $t$ approaches $0\in\Delta$. As we have all freedom of choice for the family of metrics $(\gamma_t)_{t\in\Delta}$, we will endeavour to find a special choice ensuring the uniform boundedness of the volumes.

 As every effective $(n-p)$-cycle $Z_t=\sum_jn_j(t)\, Z_j(t)$ on $X_t$ is a finite linear combination with positive integers $n_j(t)$ of irreducible analytic subsets $Z_t\subset X_t$ of dimension $n-p$, the associated De Rham cohomology class $\{[Z_t]\}\in H^{2p}(X, \, \R)$ is {\it integral}. Thus the map 

$$\Delta^{\star}\ni t\mapsto \{[Z_t]\}\in H^{2p}(X, \, \Z),$$

\noindent being continuous and integral-class-valued, must be constant. Fix any {\it real} ($d$-closed) differential $(2p)$-form $\alpha$ in this constant De Rham class. As $[Z_t]$ and $\alpha$ are $d$-cohomologous for every $t\in\Delta^{\star}$, there exists a {\it real} current $\beta'_t$ of degree $(2p-1)$ on $X$ such that

\begin{equation}\label{eqn:beta'}\alpha=[Z_t] + d\beta'_t, \hspace{3ex} t\in\Delta^{\star}.\end{equation}

 A double upper index $r, s$ will denote throughout the component of pure type $(r, \, s)$ of the form or current to which the index is attached. Since the current $[Z_t]$ is of pure type $(p, \, p)$, identifying the pure-type components on either side of the equality, we see that identity (\ref{eqn:beta'}) is equivalent to the following set of identities  for all $t\in\Delta^{\star}$: 

\begin{eqnarray}\label{eqnarray:puretype'}\nonumber\alpha_t^{0, \, 2p} & = & \bar\partial_t\beta^{'0, \, 2p-1}_t,\\
\nonumber\alpha_t^{1, \, 2p-1}-\partial_t\beta^{'0, \, 2p-1}_t=\bar\partial_t\beta^{'1, \, 2p-2}_t, & \dots & , \alpha_t^{p-1, \, p+1}-\partial_t\beta^{'p-2, \, p+1}_t=\bar\partial_t\beta^{'p-1, \, p}_t, \\
\nonumber\alpha_t^{p, \, p}-\partial_t\beta^{'p-1, \, p}_t-[Z_t] & = & \bar\partial_t\beta^{'p, \, p-1}_t,\\
\nonumber\alpha_t^{p+1, \, p-1}-\partial_t\beta^{'p, \, p-1}_t=\bar\partial_t\beta^{'p+1, \, p-2}_t, & \cdots & , \alpha_t^{2p-1, \, 1}-\partial_t\beta^{'2p-2, \, 1}_t=\bar\partial_t\beta^{'2p-1, \, 0}_t, \\
 \alpha_t^{2p, \, 0} & = & \partial_t\beta^{'2p-1,\, 0}_t.\end{eqnarray}

\noindent For all $t\in\Delta^{\star}$, we also have $\beta'_t= \overline{\beta'_t}$ (as $\beta'_t$ is real) which amounts to 

\begin{equation}\beta_t^{'l, \, 2p-1-l}=\overline{\beta_t^{'2p-1-l, \, l}}, \hspace{3ex} l=0, 1, \dots , 2p-1.\end{equation}

 The current $\beta'_t$ is determined only up to the kernel of $d$. We now proceed to construct a {\it real} $C^{\infty}$ $(2p-1)$-form $\beta_t$, having the same properties as the current $\beta'_t$, by inductively choosing its pure-type components to be minimal $L^2$-norm solutions (w.r.t. $\gamma_t$) of the first half of equations (\ref{eqnarray:puretype'}) for all $t\in\Delta^{\star}$.

 Thus, for every $t\in\Delta^{\star}$, let $\beta^{0, \, 2p-1}_t$ be the form of $J_t$-type $(0, \, 2p-1)$ which is the minimal $L^2$-norm solution of the equation (cf. first equation in (\ref{eqnarray:puretype'})):

\begin{equation}\label{eqn:beta02p-1}\alpha_t^{0, \, 2p}=\bar\partial_t\beta^{0, \, 2p-1}_t, \hspace{3ex} t\in\Delta^{\star}.\end{equation}

 In other words, $\beta^{0, \, 2p-1}_t$ corrects $\beta^{'0, \, 2p-1}_t$ if the latter is not of minimal $L^2$-norm among the solutions of the above equation. We have an explicit formula for the minimal $L^2$-norm solution:

\begin{equation}\label{eqn:minsol0}\beta^{0, \, 2p-1}_t=\Delta_t^{''-1}\bar\partial_t^{\star}\alpha_t^{0, \, 2p}, \hspace{3ex} t\in\Delta^{\star},\end{equation}

\noindent where $\Delta_t'':=\bar\partial_t\bar\partial_t^{\star} + \bar\partial_t^{\star}\bar\partial_t$ denotes the $\bar\partial_t$-Laplacian defined by the metric $\gamma_t$ (involved in the adjoints) on the fibre $X_t$ for all $t\in\Delta$, while $\Delta_t^{''-1}$ denotes the inverse of the restriction of $\Delta_t''$ to the orthogonal complement of the kernel of $\Delta_t''$ (i.e. $\Delta_t^{''-1}$ is the Green operator of $\Delta_t''$). 

 To continue, we first need to ensure that $\alpha_t^{1, \, 2p-1}-\partial_t\beta_t^{0, \, 2p-1}$ is $\bar\partial_t$-exact. Given that $\alpha_t^{1, \, 2p-1}-\partial_t\beta_t^{'0, \, 2p-1}$ is $\bar\partial_t$-exact (see the second equation in (\ref{eqnarray:puretype'})), the $\bar\partial_t$-exactness of the former form is equivalent to the $\bar\partial_t$-exactness of the difference of these two forms, i.e.  the $\bar\partial_t$-exactness of:

$$(\alpha_t^{1, \, 2p-1}-\partial_t\beta_t^{0, \, 2p-1})-(\alpha_t^{1, \, 2p-1}-\partial_t\beta_t^{'0, \, 2p-1})= \partial_t(\beta_t^{'0, \, 2p-1}-\beta_t^{0, \, 2p-1}).$$

\noindent Now $d[\partial_t(\beta_t^{'0, \, 2p-1}-\beta_t^{0, \, 2p-1})]=0$ because $\partial_t^2(\beta_t^{'0, \, 2p-1}-\beta_t^{0, \, 2p-1})=0$ and 

$$\bar\partial_t\partial_t(\beta_t^{'0, \, 2p-1}-\beta_t^{0, \, 2p-1})=-\partial_t\bar\partial_t(\beta_t^{'0, \, 2p-1}-\beta_t^{0, \, 2p-1})=-\partial_t(\alpha_t^{0, \, 2p}- \alpha_t^{0, \, 2p})=0,$$

\noindent thanks to the fact that $\bar\partial_t\beta^{'0, \, 2p-1}_t=\bar\partial_t\beta^{0, \, 2p-1}_t$ as both $\beta^{0, \, 2p-1}_t$ and $\beta^{'0, \, 2p-1}_t$ are solutions of equation (\ref{eqn:beta02p-1}) (see also the first equation in (\ref{eqnarray:puretype'})). Thus the pure type $(1, \, 2p-1)$-form $\partial_t(\beta_t^{'0, \, 2p-1}-\beta_t^{0, \, 2p-1})$ is $d$-closed and also, in an obvious way, $\partial_t$-exact for all $t\in\Delta^{\star}$. Then the $\partial\bar\partial$-lemma (which has been supposed to hold on $X_t$ for $t\neq 0$) implies the $\bar\partial_t$-exactness of $\partial_t(\beta_t^{'0, \, 2p-1}-\beta_t^{0, \, 2p-1})$ for all $t\neq 0$. This in turn implies, as has already been argued, that $\alpha_t^{1, \, 2p-1}-\partial_t\beta_t^{0, \, 2p-1}$ is $\bar\partial_t$-exact for all $t\in\Delta^{\star}$.  
 
 Considering now the analogue of the second equation in (\ref{eqnarray:puretype'}), we define $\beta_t^{1, \, 2p-2}$ to be the $(2p-1)$-form of pure $J_t$-type $(1, \, 2p-2)$ which is the minimal $L^2$-norm solution of the equation:

\begin{equation}\label{eqn:beta12p-2}\alpha_t^{1, \, 2p-1}-\partial_t\beta^{0, \, 2p-1}_t=\bar\partial_t\beta_t^{1, \, 2p-2}, \hspace{3ex} t\in\Delta^{\star}.\end{equation}

\noindent This equation does have solutions since we have proved that its left-hand side is $\bar\partial_t$-exact for all $t\in\Delta^{\star}$. We can thus go on inductively to construct forms $\beta^{l, \, 2p-1-l}_t$ of $J_t$-type $(l, \, 2p-1-l)$ for all $l\in\{0, 1, \dots , p-1\}$ and all $t\in\Delta^{\star}$. Indeed, once $\beta^{l-1, \, 2p-l}_t$ has been constructed as the minimal $L^2$-norm solution of the equation

\begin{equation}\label{eqn:betal-12p-l}\alpha_t^{l-1, \, 2p-l+1}-\partial_t\beta^{l-2, \, 2p-l+1}_t=\bar\partial_t\beta_t^{l-1, \, 2p-l}, \hspace{3ex} t\in\Delta^{\star},\end{equation} 

\noindent the pure-type form $\alpha_t^{l, \, 2p-l}-\partial_t\beta^{l-1, \, 2p-l}_t$ is seen to be $\bar\partial_t$-exact by the same argument using the $\partial\bar\partial$-lemma on $X_t$ ($t\neq 0$) as the one spelt out above for $l=1$. The form $\beta^{l, \, 2p-l-1}_t$ is then defined to be the minimal $L^2$-norm solution of the equation

\begin{equation}\label{eqn:betal+12p-2-l}\alpha_t^{l, \, 2p-l}-\partial_t\beta^{l-1, \, 2p-l}_t=\bar\partial_t\beta_t^{l, \, 2p-l-1}, \hspace{3ex} t\in\Delta^{\star}.\end{equation}

 In this case, the explicit formula giving the minimal solution reads

\begin{equation}\label{eqn:minsoll}\beta^{l, \, 2p-l-1}_t=\Delta_t^{''-1}\bar\partial_t^{\star}(\alpha_t^{l, \, 2p-l}-\partial_t\beta^{l-1, \, 2p-l}_t), \hspace{3ex} t\in\Delta^{\star}, \,\,\, l=1, \dots p-1,\end{equation}

\noindent where $\Delta_t'':C^{\infty}_{l, \, 2p-l-1}(X_t, \, \C)\rightarrow C^{\infty}_{l, \, 2p-l-1}(X_t, \, \C)$ is the $\bar\partial_t$-Laplacian defined on the space of $(l, \, 2p-l-1)$-forms of class $C^{\infty}$ on $X_t$, as recalled earlier.

 In this fashion we have defined smooth forms $\beta^{0, \, 2p-1}_t, \beta^{1, \, 2p-2}_t, \dots ,  \beta^{p-1, \, p}_t$ for all $t\in\Delta^{\star}$. They satisfy the first $p$ equations (with $\beta_t$ replacing $\beta_t'$) among the $(2p+1)$ equations in (\ref{eqnarray:puretype'}). We then go on to define, for all $t\in\Delta^{\star}$, smooth forms $\beta^{p, \, p-1}_t, \beta^{p+1, \, p-2}_t, \dots , \beta^{2p-1, \, 0}_t$ as the conjugates of the previous set of forms taken in reverse order:

\begin{equation}\label{eqn:betaconj}\beta^{p+s, \, p-s-1}_t:=\overline{\beta^{p-s-1, \, p+s}_t}, \hspace{3ex} s=0, 1, \dots , p-1, \,\, t\in\Delta^{\star}.\end{equation} 

 Since the form $\alpha$ has been chosen to be real, we take conjugates and see that the forms $\beta^{p+s, \, p-s-1}_t$, $s=0, 1, \dots , p-1,$ satisfy the last $p$ equations (with $\beta_t$ replacing $\beta_t'$) among the $(2p+1)$ equations in (\ref{eqnarray:puretype'}). If we now set

\begin{equation}\label{eqn:beta}\beta_t:=\beta^{0, \, 2p-1}_t + \cdots + \beta^{p-1, \, p}_t + \beta^{p, \, p-1}_t + \cdots + \beta^{2p-1, \, 0}_t, \hspace{3ex} t\in\Delta^{\star},\end{equation} 

\noindent we obtain a family $(\beta_t)_{t\in\Delta^{\star}}$ of real $C^{\infty}$ forms of degree $2p-1$ on $X$ varying in a $C^{\infty}$ way with $t\in\Delta^{\star}$. Moreover, the $(2p)$-current $\alpha-[Z_t]-d\beta_t$ is of pure type $(p, \, p)$ for all $t\in\Delta^{\star}$ as can be seen from the construction of $\beta_t$: its pure-type components satisfy the analogues for $\beta_t$ (instead of $\beta_t'$) of equations (\ref{eqnarray:puretype'}), except the one involving $[Z_t]$, which amount to the vanishing of all the pure-type components of $\alpha-[Z_t]-d\beta_t$, except the one of type $(p, \, p)$ which is the only one to which $[Z_t]$ contributes. The current $\alpha-[Z_t]-d\beta_t$ is also $d$-exact in an obvious way (it equals $d(\beta_t'-\beta_t)$).

 A final application of the $\partial\bar\partial$-lemma (supposed to hold on $X_t$ for every $t\neq 0$) shows that $\alpha-[Z_t]-d\beta_t$ is also $\partial_t\bar\partial_t$-exact for $t\neq 0$. Thus there exists a family $(R_t)_{t\in\Delta^{\star}}$ of $(2p-2)$-currents of respective $J_t$-types $(p-1, \, p-1)$ such that

\begin{equation}\label{eqn:ddbarcurrent}\alpha=[Z_t] + d\beta_t + \partial_t\bar\partial_t R_t, \hspace{3ex} t\in\Delta^{\star}.\end{equation}

\begin{Conc}\label{Conc:vol} If the $\partial\bar\partial$-lemma holds on $X_t$ for all $t\in\Delta^{\star}$, the $\gamma_t$-volumes of any $C^{\infty}$ family of relative $(n-p)$-cycles $(Z_t)_{t\in\Delta^{\star}}$ can be expressed as

\begin{equation}\label{eqn:volexp}v_{\gamma_t}(Z_t):=\int\limits_X[Z_t]\wedge\gamma_t^{n-p}= \int\limits_X\alpha\wedge\gamma_t^{n-p} - \int\limits_Xd\beta_t\wedge\gamma_t^{n-p} -  \int\limits_X\partial_t\bar\partial_t R_t\wedge\gamma_t^{n-p},  \hspace{3ex} t\in\Delta^{\star},\end{equation}

\noindent for any family of Hermitian metrics $(\gamma_t)_{t\in\Delta}$ on the fibres $(X_t)_{t\in\Delta}$, where $\alpha$ is a fixed real $(2p)$-form in the De Rham class that is common to all $[Z_t]$, $(\beta_t)_{t\in\Delta^{\star}}$ are given by formula (\ref{eqn:beta}) by adding their components inductively defined in formulae (\ref{eqn:minsol0}), (\ref{eqn:minsoll}) and (\ref{eqn:betaconj}), while $(R_t)_{t\in\Delta^{\star}}$ are given by (\ref{eqn:ddbarcurrent}).

\end{Conc}

 Recall that what is at stake is ensuring that $v_{\gamma_t}(Z_t)$ is uniformly bounded as $t\in\Delta^{\star}$ approaches $0\in\Delta$. If $\gamma_t$ is chosen to vary in a $C^{\infty}$ way with $t\in\Delta$ (up to $t=0$), the first term in the right-hand side of (\ref{eqn:volexp}) stays bounded when $t$ varies in a relatively compact neighbourhood $U\Subset\Delta$ of $0\in\Delta$, since $\alpha$ is independent of $t$. The other two terms are problematic as both $\beta_t$ and $R_t$ are only defined off $t=0\in\Delta$.   

 The first observation is that, when the cycles $Z_t$ are divisors (i.e. $p=1$), the third term in the right-hand side of (\ref{eqn:volexp}) can be easily handled. The reason is that the Hermitian metrics $\gamma_t$ of the fibres $X_t$ can be chosen as {\it Gauduchon metrics}, i.e. such that $\partial_t\bar\partial_t\gamma_t^{n-1}=0$ for all $t\in\Delta$. Indeed, Gauduchon metrics exist on every compact complex manifold (cf. [Gau77]) and, moreover, one can always find a family $(\gamma_t)_{t\in\Delta}$, varying in a $C^{\infty}$ way with $t$, of Gauduchon metrics on the fibres $(X_t)_{t\in\Delta}$ of any smooth holomorphic family of compact complex manifolds. The argument for this last (well-known) assertion is recalled, for instance, in [Pop09, section 2]. With this special choice for $(\gamma_t)_{t\in\Delta}$, Stokes' theorem gives:

$$\int\limits_X\partial_t\bar\partial_t R_t\wedge\gamma_t^{n-1}=-\int\limits_X R_t\wedge\partial_t\bar\partial_t\gamma_t^{n-1}=0, \hspace{3ex} t\in\Delta,$$

\noindent so this term vanishes in the case of divisors. However, achieving uniform boundedness for this term in the case of higher codimensional cycles (i.e. for $p\geq 2$) is a major challenge.

As for uniformly bounding the term depending on $\beta_t$ in the right-hand side of (\ref{eqn:volexp}), the difficulty stems from the possible jump of the Hodge numbers $h^{p, \, q}(t):=\mbox{dim}_{\C}H^{p, \, q}(X_t, \, \C)$ at $t=0$. The family of strongly elliptic operators $(\Delta_t'')_{t\in\Delta}$ defined in $J_t$-bidegree $(p, \, q)$ varies in a $C^{\infty}$ way with $t$, while a classical result of Kodaira and Spencer [KS60] ensures that the corresponding family of Green operators $(\Delta_t^{''-1})_{t\in\Delta}$ varies in a $C^{\infty}$ with $t$ if the dimension (as a $\C$-vector space) of the kernel $\ker\Delta_t''$ is independent of $t\in\Delta$. Since $\ker\Delta_t''$ is isomorphic to the Dolbeault cohomology space $H^{p, \, q}(X_t, \, \C)$ by the Hodge Isomorphism Theorem, we have differentiability of the families of operators $(\Delta_t^{''-1})_{t\in\Delta}$ (and hence of the families of forms $(\beta_t^{l, \, 2p-l-1})_{t\in\Delta}$, $l=, 1, \dots , p-1$, thanks to the formulae (\ref{eqn:minsol0}) and (\ref{eqn:minsoll})) if the Hodge numbers $h^{l, \, 2p-l-1}(t)$, $l=, 1, \dots , p-1$ , of the fibres do not jump at $t=0\in\Delta$. This condition is fulfilled, for instance, under the hypothesis of Conjecture \ref{Conj:relpropC} since the {\it class} ${\cal C}$ assumption on the fibres ensures the degeneracy at $E_1^{\bullet}$ of the Fr\"olicher spectral sequence of each fibre which, in turn, is known to imply local constancy of the Hodge numbers of the fibres. Thus the term depending on $\beta_t$ in the expression (\ref{eqn:volexp}) for $v_{\gamma_t}(Z_t)$ is uniformly bounded when $t$ varies in a relatively compact neighbourhood $U\Subset\Delta$ of $0\in\Delta$ under the hypothesis of Conjecture \ref{Conj:relpropC}. However, controlling this term in the more general situation of Conjecture \ref{Conj:relpropC*} poses a major challenge as the Hodge numbers might {\it a priori} jump at $t=0$ if the {\it class} ${\cal C}$ assumption skips $X_0$ (unless they can be shown not to do so, which seems a daunting task).
    
 A by-product of these considerations is that the divisor case of Conjecture \ref{Conj:relpropC} holds true.

\vspace{2ex}

\noindent {\it End of proof of Proposition \ref{Prop:reldivprop}.} The statement of Proposition \ref{Prop:reldivprop} falls into the mould of Conjecture \ref{Conj:relpropC} (even with the weaker $\partial\bar\partial$-lemma assumption over the fibres above $\Delta^{\star}$ replacing the {\it class} ${\cal C}$ one), but only deals with the special case of divisors. Thus $p=1$ and the term depending on $R_t$ in the expression (\ref{eqn:volexp}) for $v_{\gamma_t}(Z_t)$ vanishes if the metrics $(\gamma_t)_{t\in\Delta}$ are chosen to be {\it Gauduchon metrics}, as has been explained above. The major challenge posed by other term, depending on $\beta_t$, in the right-hand side of (\ref{eqn:volexp}) has been solved in [Pop09]. Before briefly recalling the argument, we wish to emphasise that the case where $p\geq 2$ falls completely outside the method of [Pop09] and of the present paper and is thus widely open.

 As $p=1$, formula (\ref{eqn:minsol0}) defining $\beta_t^{0, \, 1}$ reads

\begin{equation}\label{eqn:minsol0bis}\beta^{0, \, 1}_t=\Delta_t^{''-1}\bar\partial_t^{\star}\alpha_t^{0, \, 2}, \hspace{3ex} t\in\Delta^{\star},\end{equation}

\noindent while $\beta_t=\overline{\beta_t^{0, \, 1}} + \beta_t^{0, \, 1}$ (cf. (\ref{eqn:betaconj}) and (\ref{eqn:beta})) is now a $1$-form. Thus only the $(1, \, 1)$-component of $d\beta_t$ has a non-trivial contribution to $v_{\gamma_t}(Z_t)$ and we get

$$\int\limits_Xd\beta_t\wedge\gamma_t^{n-1}=\int\limits_X(\partial_t\beta_t^{0, \, 1} + \bar\partial_t\beta_t^{1, \, 0})\wedge\gamma_t^{n-1},$$

\noindent where we have set $\beta_t^{1, \, 0}:=\overline{\beta_t^{0, \, 1}}$. As $\partial_t\beta_t^{0, \, 1}$ and $\bar\partial_t\beta_t^{1, \, 0}$ are conjugate to each other, it suffices to uniformly bound

\begin{equation}\label{eqn:mainq}I_t:=\int\limits_X\partial_t\beta_t^{0, \, 1}\wedge\gamma_t^{n-1}, \hspace{3ex} t\in\Delta^{\star}.\end{equation}

 The difficulty is that $\beta_t^{0, \, 1}$ (hence also $\partial_t\beta_t^{0, \, 1}$) might {\it explode} as $t\in\Delta^{\star}$ approaches $0\in\Delta$, if $h^{0, \, 1}(t)$ {\it jumps} at $t=0$. However, $\bar\partial_t\beta_t^{0, \, 1}=\alpha_t^{0, \, 2}$ (see equation (\ref{eqn:beta02p-1}) with $p=1$) and thus $\bar\partial_t\beta_t^{0, \, 1}$ extends in a $C^{\infty}$ way to $t=0$ since the $(0, \, 2)$-component $\alpha_t^{0, \, 2}$ of the fixed form $\alpha$ w.r.t. to the holomorphic family of complex structures $(J_t)_{t\in\Delta}$ does. Hence the idea of trying to substitute $\bar\partial_t\beta_t^{0, \, 1}$ for $\partial_t\beta_t^{0, \, 1}$ in (\ref{eqn:mainq}) appears as natural. Stokes' theorem gives

\begin{equation}\label{eqn:mainqS}I_t=\int\limits_X\beta_t^{0, \, 1}\wedge\partial_t\gamma_t^{n-1}, \hspace{3ex} t\in\Delta^{\star}.\end{equation}

 Recall that the metrics $\gamma_t$, $t\in\Delta$, have been chosen to satisfy the {\it Gauduchon condition}: $\partial_t\bar\partial_t\gamma_t^{n-1}=0$ for all $t\in\Delta$. Thus $d(\partial_t\gamma_t^{n-1})=0$, $t\in\Delta$, and the $\partial\bar\partial$-lemma (supposed to hold on every $X_t$ with $t\neq 0$) implies that the $d$-closed form $\partial_t\gamma_t^{n-1}$ of pure type $(n, \, n-1)$, which is obviously $\partial_t$-exact, must also be $\bar\partial_t$-exact for every $t\neq 0$. However, it is not clear {\it a priori} whether $\partial_0\gamma_0^{n-1}$ is $\bar\partial_0$-exact since the $\partial\bar\partial$-lemma is not known to hold on $X_0$. According to [Pop09, Definition 3.1], a Hermitian metric $\gamma_0$ on $X_0$ is said to be a {\it strongly Gauduchon metric} if 

\begin{equation}\label{eqn:sGdef}\partial_0\gamma_0^{n-1} \hspace{2ex} \mbox{is}\hspace{2ex} \bar\partial_0-\mbox{exact},\end{equation}

\noindent i.e. if the condition that is needed here is met. (This condition clearly implies the {\it Gauduchon condition}). We have shown in [Pop09, Lemma 3.2] that, although a {\it strongly Gauduchon metric} need not exist on an arbitrary compact complex manifold, the existence of such a metric $\gamma_0$ on $X_0$ is equivalent to the existence of a real $d$-closed $C^{\infty}$ form $\Omega$ of degree $2n-2$ on $X$ such that its component of $J_0$-type $(n-1, \, n-1)$ is positive definite (i.e. $\Omega_0^{n-1, \, n-1}>0$). Now, if $X_0$ carries a {\it strongly Gauduchon metric} $\gamma_0$, the components $\Omega_t^{n-1, \, n-1}$ of $J_t$-type $(n-1, \, n-1)$ of $\Omega$ vary in a $C^{\infty}$ way with $t\in\Delta$ and, therefore, the strict positivity condition is preserved in a small neighbourhood of $0\in\Delta$ (and thus on the whole $\Delta$ if $\Delta$ is shrunk sufficiently about $0$):

$$\Omega_t^{n-1, \, n-1}>0, \hspace{3ex} t\in\Delta.$$

\noindent Thus $\Omega$ defines a {\it strongly Gauduchon metric} on every fibre $X_t$ with $t\in\Delta$ (after possibly shrinking $\Delta$). This shows that the {\it strongly Gauduchon condition} is open in the classical topology of the base under holomorphic deformations. Moreover, since the form $\Omega$ is {\it real}, the closedness condition $d\Omega=0$ is equivalent to

$$\partial_t\Omega_t^{n-1, \, n-1}=-\bar\partial_t\Omega_t^{n, \, n-2}, \hspace{3ex} t\in\Delta.$$

\noindent Thus the $\bar\partial_t$-potentials $\Omega_t^{n, \, n-2}$ of $\partial_t\Omega_t^{n-1, \, n-1}$ also vary in a $C^{\infty}$ way with $t\in\Delta$ since they are components of pure $J_t$-type $(n, \, n-2)$ of the fixed form $\Omega$. 

\begin{Conc}\label{Conc:smoothfampot} Let $\pi:{\cal X}\rightarrow\Delta$ be an arbitrary complex analytic family of compact complex manifolds. Suppose that $X_0$ carries a {\bf strongly Gauduchon metric} $\gamma_0$. Then, after possibly shrinking $\Delta$ about $0$, there exists a family $(\gamma_t)_{t\in\Delta}$, varying in a $C^{\infty}$ way with $t$, of {\bf strongly Gauduchon metrics} on the respective fibres $(X_t)_{t\in\Delta}$. Moreover, there exists a family $(\zeta_t^{n, \, n-2})_{t\in\Delta}$, varying in a $C^{\infty}$ way with $t$, of $(2n-2)$-forms on $X$ of respective $J_t$-types $(n, \, n-2)$ such that

$$\partial_t\gamma_t^{n-1}=\bar\partial_t\zeta_t^{n, \, n-2}, \hspace{3ex} t\in\Delta.  $$

\end{Conc}

 Here the emphasis is on the differentiable dependence of $\zeta_t^{n, \, n-2}$ on $t\in\Delta$. Clearly, the link between the $(2n-2)$-form $\Omega$ mentioned above and the objects of Conclusion \ref{Conc:smoothfampot} is

$$\gamma_t^{n-1}=\Omega_t^{n-1, \, n-1}, \hspace{3ex} \zeta_t=-\Omega_t^{n, \, n-2}, \hspace{3ex} t\in\Delta.$$

 To conclude, we now need the following crucial ingredient from [Pop09].

\begin{Prop}\label{Prop:re} (Proposition 4.1 in [Pop09]) If the $\partial\bar\partial$-lemma holds on $X_t$ for every $t\in\Delta^{\star}$, then $X_0$ carries a {\bf strongly Gauduchon metric}. 

\end{Prop}

 This means that under the hypothesis of Proposition \ref{Prop:reldivprop}, Conclusion \ref{Conc:smoothfampot} holds and, choosing a differentiable family $(\gamma_t)_{t\in\Delta}$ of {\it strongly Gauduchon metrics} on the fibres $(X_t)_{t\in\Delta}$, (\ref{eqn:mainqS}) reads

\begin{eqnarray}\label{eqn:finalest}I_t & = & \int\limits_X\beta_t^{0, \, 1}\wedge\partial_t\gamma_t^{n-1}=\int\limits_X\beta_t^{0, \, 1}\wedge\bar\partial_t\zeta_t^{n, \, n-2}, \\
 & = & \int\limits_X\bar\partial_t\beta_t^{0, \, 1}\wedge\zeta_t^{n, \, n-2}=\int\limits_X\alpha_t^{0, \, 2}\wedge\zeta_t^{n, \, n-2}, \hspace{3ex} t\in\Delta^{\star},\end{eqnarray}

\noindent where Stokes' theorem has been applied in passing to the second line. As both families of forms $(\alpha_t^{0, \, 2})_{t\in\Delta}$ and $(\zeta_t^{n, \, n-2})_{t\in\Delta}$ vary in a $C^{\infty}$ way with $t$ (up to $t=0$), $I_t$ is bounded independently of $t\in\Delta^{\star}$ after possibly shrinking $\Delta$ about $0$. Hence the volume $v_{\gamma_t}(Z_t)$ is bounded independently of $t$ when $t\in\Delta^{\star}$ approaches $0\in\Delta$ (see (\ref{eqn:volexp})). 

\vspace{2ex}

 To show properness over $\Delta$ of an arbitrary irreducible component $S\subset{\cal C}^{n-1}({\cal X}/\Delta)$, one has to show that for every compact subset $K\subset\Delta$, $\mu_{n-1}^{-1}(K)\cap S$ is a compact subset of ${\cal C}^{n-1}({\cal X}/\Delta)$. If $(Z_s)_{s\in S}$ is the analytic family of divisors associated with $S$ (such that $Z_s\subset X_{\mu_{n-1}(s)}$, $s\in S$), this amounts to proving that the volumes

$$v_{\gamma_s}(Z_s)=\int\limits_X[Z_s]\wedge\gamma_s^{n-1}$$

\noindent are uniformly bounded when $s$ ranges over $\mu_{n-1}^{-1}(K)\cap S$. Here we have denoted for convenience $\gamma_s=\gamma_{\mu_{n-1}(s)}$. As mentioned in the Introduction, the absolute Barlet space ${\cal C}^{n-1}(X_t)$ of divisors of every fibre $X_t$ is known to have compact irreducible components. Thus $
v_{\gamma_s}(Z_s)$ stays uniformly bounded when $Z_s$ varies across any irreducible component of any given fibre. It then suffices to show uniform boundedness of the volumes in the {\it horizontal directions}, i.e. when $Z_t\subset X_t$ varies in a differentiable family $(Z_t)_{t\in\Delta^{\star}}$ with $t\in\Delta^{\star}$ approaching $0\in\Delta$. This has been done above. The proof of Proposition \ref{Prop:reldivprop} is complete. \hfill $\Box$

\vspace{6ex}

\noindent {\bf References.} \\

\noindent [Bar75]\, D. Barlet --- {\it Espace analytique r\'eduit des cycles analytiques complexes compacts d'un espace analytique complexe de dimension finie} --- Fonctions de plusieurs variables complexes, II (S\'em. Fran\c{c}ois Norguet, 1974-1975), LNM, Vol. {\bf 482}, Springer, Berlin (1975) 1-158.

\vspace{1ex}

\noindent [Bis64]\, E. Bishop --- {\it Conditions for the Analyticity of Certain Sets} --- Mich. Math. J. {\bf 11} (1964) 289-304.

\vspace{1ex}

\noindent [Cam80]\, F. Campana --- {\it Alg\'ebricit\'e et compacit\'e dans l'espace des cycles d'un espace analytique complexe} --- Math. Ann. {\bf 251} (1980), 7-18.

\vspace{1ex}

\noindent [CP94]\, F. Campana, T. Peternell --- {\it Cycle spaces} --- Several Complex Variables, VII,  319-349, Encyclopaedia Math. Sci., {\bf 74}, Springer, Berlin (1994).

\vspace{1ex}

\noindent [Dem 97] \, J.-P. Demailly --- {\it Complex Analytic and Algebraic Geometry}---http://www-fourier.ujf-grenoble.fr/~demailly/books.html

\vspace{1ex}

\noindent [Fuj78]\, A. Fujiki --- {\it Closedness of the Douady Spaces of Compact K\"ahler Spaces} --- Publ. RIMS, Kyoto Univ. {\bf 14} (1978), 1-52.

\vspace{1ex}

\noindent [FP09]\, A. Fujiki, M. Pontecorvo --- {\it Non-Upper-Semicontinuity of Algebraic Dimension for Families of Compact Complex Manifolds} --- arXiv e-print math.AG/0903.4232v2.

\vspace{1ex}

\noindent [Gau77]\, P. Gauduchon --- {\it Le th\'eor\`eme de l'excentricit\'e nulle} --- C.R. Acad. Sc. Paris, S\'erie A, t. {\bf 285} (1977), 387-390.



\vspace{1ex}

\noindent [Kod86]\, K. Kodaira --- {\it Complex Manifolds and Deformations of Complex Structures} --- Grundlehren der Math. Wiss. {\bf 283}, Springer (1986).

\vspace{1ex}

\noindent [KS60]\, K. Kodaira, D.C. Spencer --- {\it On Deformations of Complex Analytic Structures. III. Stability Theorems for Complex Structures} --- Ann. of Math. (2) {\bf 71} No.1 (1960) 43-76.

\vspace{1ex}

\noindent [Lie78]\, D. Lieberman --- {\it Compactness of the Chow Scheme: Applications to Automorphisms and Deformations of K\"ahler Manifolds} --- Lect. Notes Math. {\bf 670} (1978), 140-186.\vspace{1ex}

\vspace{1ex}

\noindent [Moi67]\, B.G. Moishezon --- {\it On $n$-dimensional Compact Varieties with $n$ Algebraically Independent Meromorphic Functions} --- Amer. Math. Soc. Translations {\bf 63} (1967) 51-177.

\vspace{1ex}

\noindent [Pop09]\, D. Popovici --- {\it Limits of Projective Manifolds Under Holomorphic Deformations} --- arXiv e-print math.AG/0910.2032v1.

\vspace{1ex}

 \noindent [Var86]\, J. Varouchas --- {\it Sur l'image d'une vari\'et\'e k\"ahl\'erienne compacte} --- LNM {\bf 1188}, Springer (1986) 245--259.

\vspace{2ex}

\noindent Universit\'e Paul Sabatier, Institut de Math\'ematiques de Toulouse, 118 route de Narbonne, 31062 Toulouse Cedex 9, France

\noindent Email: popovici@math.ups-tlse.fr

\end{document}